\documentclass[10pt,a4paper,twoside]{article}
\usepackage{amssymb}
\usepackage{latexsym}
\usepackage{fullpage}
\usepackage{a4,psfig,latexsym,rotating,amsmath,epsfig}
\usepackage{here}
\usepackage{amsmath}
\usepackage{amscd}

\def\R{\mathbb{R}}

\def\D{\mathbb{D}}

\def\E{\mathfrak{E}}

\def\1{\mathbf{1}}
\def\:{\lrcorner}
\def\#{\sharp}

\def\o{\circ}

\def\<#1,#2>{\langle#1,\,#2\rangle}
\def\E{\rm{Emb}}

\def\S{\mathbb{S}\,}

\def\qed{\ensuremath{\quad\Box\quad}}
\def\pfill{\par\vskip2mm plus1mm minus1mm\noindent}
\def\inv#1{\raise.1em\hbox to 0pt{$^{-1}$\hss}_{#1}\;}

\newtheorem{Theorem}{Theorem}
\newtheorem{Lemma}[Theorem]{Lemma}
\newtheorem{Corollary}[Theorem]{Corollary}

\begin{document}

\title{A note on $H^1 ( \E (M,N) )$}

\author{Olaf M\"uller}
\maketitle
\begin{center}Max-Planck-Institute for Mathematics in the Sciences\\Inselstrasse 22-26\\D-04103 Leipzig, Europe\end{center}
\date{}

\begin{abstract}
\noindent Aim of this note is to gain cohomological information about the manifold $ \E (M,N)$ from the topology of the target manifold $N$. For special conditions, a monomorphism $H^1(N) \rightarrow H^1(Emb(M,N))$ is constructed.
\end{abstract}

\noindent Let $\E_g (M,N)$ be the manifold of all embeddings of a manifold $M$ into a manifold $N$ that differ from a fixed reference embedding $g$ only within a compact subset of $M$. Note that if $M$ is compact then $\E_g (M,N) = \E (M,N)$, the space of all embeddings. Now and in the following $M$ should be nonempty and $N$ at least one-dimensional to give sense to the notion of 1-forms. $\E_g (M,N)$ is infinite-dimensional iff $dim(M) \geq 1$.

\bigskip

\noindent {\em Note.} The following theorem \ref{zentral} is a straightforward consequence within the framework of fiber integration (cf. ~\cite{sw}, ~\cite{abr}), nevertheless it may be instructive to see another proof for the particular case.

\bigskip

\noindent For our construction, a crucial tool is the identification of a tangent vector $V$ resp. the value of a vector field $V$ on $ \E_g (M,N)$ at a fixed map $f$ with a vector field {\em along} $f$, i.e. a section of $f^*TN$:

$$\hat{} \  \vert_{f} :{\rm Vect} (\E_g (M,N)) \rightarrow \Gamma (f^*TN),$$

$$\hat{V} \vert_{f} (p) := \mathcal{L}_V ev_p$$

or equivalently,

$$ \hat{V} \vert_{f} : p \mapsto \partial_t (f(p)),$$

\noindent where $p \in M$, $f_t$ a curve representing $V( f)$. This means, we just fix a point $p \in M$ and note the direction in which it is moved infinitesimally by the family of maps $f_t$. $\hat{V} \vert_{f}$ has compact support because $ \bigcup_{t \in [-1;1]} \{ f_t \neq g \}$ is compact.

\bigskip

\noindent Conversely, a vector field $V$ on $N$ maps to a vector field $ \tilde{V}$ on $ \E_g (M,N)$ by 

$$ \tilde{} : {\rm Vect} (N)  \rightarrow {\rm Vect} ( \E_g (M,N))$$

$$ \tilde{V}: f \mapsto \partial_t \vert_{t=0} (Fl_t^V \o f) ,$$

where $Fl_t^V $ is the flow of $V$.

\begin{Lemma}

\label{correspondence}

\begin{enumerate}

\item[(i)]
Both maps $\ \tilde{} \ $ and $\ \hat{}\ \vert_{f}  $ are smooth for all $f$; 

moreover, $\  \hat{} \ $ depends smoothly on $ f $.

\item[(ii)]
$ \ \hat{} \  \vert_{f} \o \ \tilde{} \   = f^*  \qquad \forall f \in \E_g (M,N) $, 

i.e. the composition is the pull-back along $f$.

\item[(iii)]
The Lie bracket is preserved by $\ \tilde{\ } \ $:

$$ [\tilde{V}, \tilde{W}] = \widetilde{[V,W]} .$$

\end{enumerate}

\end{Lemma}

\noindent {\em Proof.} The infinite-dimensional manifold of vector fields on a manifold is in particular a Fr\^echet space, and the flow of a vector field depends smoothly on it, so both maps are smooth. Looking deeply at the definitions, one sees that they are inverse to each other in the sense of (ii). For a proof of the third statement we define

$$ A: {\rm Diff}(N) \rightarrow {\rm Diff}( \E_g (M,N)) ,$$

$$ h \mapsto ( f \mapsto h \o f) .$$

\noindent As this map is a differentiable Lie group homomorphism, its differential  
$dA \vert_{\1} : T_{\1} {\rm Diff}(N) \rightarrow T_{\1}{\rm Diff}( \E_g (M,N)) $ preserves the Lie bracket in the sense that 

\noindent $  [dA (V), dA(W)] = dA [V,W] $. Now the source resp. target of $dA$ is as Lie algebra isomorphic to ${\rm Vect}(N)$ resp. ${\rm Vect}(\E_g (M,N))$ via $\Phi$ resp. $\Psi$, and $ \Psi \o dA \o \Phi^{-1} = \  \tilde{} \ $ by the chain rule     \qed

\bigskip

\noindent In order to get a p-form $\hat{\omega}$ on $ \E_g (M,N)$ from one on $N$, for each $p$ vector fields $X_1,..X_p \in {\rm Vect} (\E_g (M,N))$ we define

$$ \hat{\omega} \vert_f (X_1,...X_p) := \int_M (\omega \o f)( \hat{X_1}\vert_{\gamma}, ...\hat{X_p}\vert_{\gamma}). $$

\noindent Note that for this definition we need a volume form on $M$.

\begin{Theorem}
\label{zentral}

$\hat{} : \Lambda^* (N) \rightarrow \Lambda^* (\E (M,N)) $  is an injective chain map. 

\end{Theorem}

\bigskip

\noindent {\em Proof.} Obviously, $\hat{\omega}$ is a p-form.
We have to show that $ d \hat{\omega} = \widehat{d \omega}$. We calculate  $ d \hat{\omega} (X_1, ... X_{p+1})$ at a  section $ f $ using

\bigskip

$ d \Omega (X_1, ... X_{n+1}) = \sum_{i=1}^{n+1}( (-1)^{i=1} X_i ( \Omega (X_1, ...X_{i-1}, X_{i+1} ... X_p))$

$   - \sum_{i<j} \Omega (-1)^{i+j} ([X_i,X_j], X_1...X_{i-1},X_{i+1}...X_{j-1}, X_{j+1}...X_n ) .$

\bigskip

\noindent Now the crucial point is that because of tensoriality of $ d \hat{\omega}$ we can choose    $X_i = \tilde{x_i} $ for some vector fields $x_i$ on $N$ . Lemma \ref{correspondence} then implies $ \widehat{[X_i,X_j]} = f^* ( [ x_i , x_j ]) $, and because of $ \hat{X_i} \vert_{f} =  f^* x_i$ we have

$$ X( \hat{\omega} (Y,Z)) \vert_f = X \int_M (\omega \o f ) ( \ \hat{Y} \vert_f , \ \hat{Z} \vert_f )        = X \int_M (\omega \o f) ( y , z) =  \int_M  x ( \omega ( y, z)) .$$

\noindent Using these equations we get the result transferring all terms to Lie brackets and derivatives in $N$.

\bigskip

\noindent Up to this point, we could have used {\em any} linear function on $ C^{\infty} (M)$ instead of the integral. Then, if $\omega$ is nonzero, there is a point $ m \in M$ and vector fields $\hat{Y}_i$ such that $ \omega \vert_{m} ( \hat{Y}_1 ... \hat{Y}_n) \neq 0$. If we scale e.g. $\hat{Y}_1$ with a function $F$ on $M$, then $\Omega ( F \hat{Y}_1 ... \hat{Y}_n ) = \int_M F \cdot  (\omega \o f) ( \hat{Y}_1 ... \hat{Y}_n)$, and the statement follows from the nondegeneracy of the scalar product of $L_2 (M)$   $\qed$ 

\bigskip

\noindent Therefore $\hat{} \ $ induces a homomorphism $\hat{} \ ^*$ of the corresponding cohomology groups. If we knew that for each exact $\hat{\omega} = d \alpha$ also the form $\omega$ is exact or, equivalently, that $\alpha = \hat{\beta}$ for some $\beta$ we would know that $\hat{} \  ^*$ is injective on the cohomology level. This is asserted for the {\em first} cohomology group by the following theorem:

\begin{Theorem}
\label{oh}
If $M$ is compact and if there is an embedding $M \rightarrow \R^{dimN}$, then
the map  $\hat{} \ ^*$  is a monomorphism $ H^1 (N) \rightarrow H^1 (\E (M,N))$.

\end{Theorem}

\noindent {\em Note.} Cohomology means here and in the following always de Rham cohomology.

\bigskip

\noindent {\em Proof.} Given a 1-form $\theta$ on $N$ we have to show that if $\hat{\theta}$ is the differential of a function $F$ on $\E (M,N) $ then $\theta$ is the differential of some function $f$ on $N$ (this implies $F = \hat{f}$ by Theorem \ref{zentral}). It is enough to show that $\int_c \theta = 0 $ for all closed curves $c$ in $N$. Given such a curve, we can construct a tubular neighborhood $T$ of $c$ and a smooth vector field on $N$ such that $Fl_V^t (T) \subset T \  \forall t$ and $c$ is the integral curve for $V$ with $Fl_V^1$ being the identity on $T$ (just take the standard circular flow vector field on the solid torus diffeomorphic to a small tubular neighborhood and scale it down radially to zero in its complement in a larger tubular neighborhood). Then the flow of $V$ on $c$ together with the distinction of a starting point $c(0)$ fixes a parametrisation of the curve to which we will refer in the following. Consider curves $K_{\delta}$ in $\rm{Emb} (M,N)$ defined by $K_{\delta}(t) = Fl_V^t \o \epsilon _{\delta} $ with $Im (\epsilon _{\delta} ) \subset B_{\delta} (c(0)) \subset T$. Then we observe

$$ \partial_t  \vert_{t_0} \int_0^t dt' \ \hat{\theta} (\dot{K}_{\delta} (t') ) = \hat{\theta} \vert_{K_{\delta} (t_0)} ( \tilde{V}) = \int_M (\theta \o K_{\delta} (t_0) ) (V \o K_{\delta} (t_0) ) = \int _M (\theta (V) ) \o K_{\delta} (t_0)$$

\noindent Therefore 

 $$ \partial_t  \vert_{t_0} \int_0^t dt' \ \hat{\theta} (\dot{K}_{\delta} (t') ) 
\underset{\delta \rightarrow 0}{\longrightarrow} vol (M) \cdot \theta (\dot{c}(t_0)) = vol(M) \cdot \partial_t  \vert_{t_0} \int_0^t dt' \  \theta ( \dot{c} (t') ) $$ 

\noindent This convergence is uniform in $t_0$: Define $D_{\delta} :=  max_{t \in [ 0 , 1 ] } diam ( Im (K_{\delta}(t)))$,   

\noindent $H :=  max_{p \in T} \parallel grad ( \theta (V) ) (p) \parallel $. Then we have 

$$\vert \int _M (\theta (V) ) \o K_{\delta} (t_0) -  vol (M) \cdot \theta (\dot{c}(t_0)) \vert \leq D_{\delta} \cdot H$$

\noindent with $D_{\delta} \rightarrow 0$ for $\delta \rightarrow 0$ which implies uniform convergence. Therefore 

\noindent $ \int_{K_{\delta}} \hat{\theta} \rightarrow vol (M) \cdot \int_c \theta$ and   $\int_c \theta = 0$  \qed

\bigskip

\noindent One easily sees that one cannot just drop the condition of the existence of an embedding $M \rightarrow \R^{dim N}$ by recalling the following fact about surfaces:

\bigskip

\noindent Let $M$ be an orientable compact two-dimensional manifold without boundary distinct from $\S^1 \times \S^1$ and $\S^2$. Then the identity component (and therefore all connected components) of the diffeomorphism group $\rm{Diff}(M)$ ($= \E (M,M)$) is homotopy equivalent to a point (cf. ~\cite{ee} , ~\cite{g} ). 

\bigskip

\noindent From Lemma 2.1. in ~\cite{sw} it is clear that $\hat{\omega}$ is invariant under the natural action of volume-preserving diffeomorphisms of $M$ on $\E (M,N)$. In the case of $M$ being the distant union of finitely many points, all diffeomorphisms (permutations) are volume-preserving if we give the same weight to each point. Therefore $\omega$ is well-defined under the usual $\Sigma_n$-action in the definition of the configuration space and we get:

\begin{Corollary} 

Let $C_m(N)$ be the n-point configuration space of a manifold $N$, $dim N \geq 1$. Then there is a monomorphism $H^1(N ) \rightarrow H^1 (C_m (N))$ \qed

\end{Corollary}

\bigskip

\noindent This result can also be compared to computations of the homology of configuration spaces as in ~\cite{bo}.

\bigskip

\noindent A result analogous to Theorem \ref{oh} for the space of differentiable maps from $M$ to $N$ (without any injectivity condition) is found in ~\cite{hv} .

\end{document}